\documentclass[10pt]{article}

\newcommand{\ddate}{\it Dedicated to Hans Duistermaat\\ on the occasion of his $65^{\mathrm th}$ birthday.}
\date{\ddate}

\usepackage{amssymb}
\usepackage[matrix,arrow,curve]{xy}
\usepackage{epsfig}

\usepackage{amsmath}

\usepackage{theorem}

\newtheorem{thm}{Theorem}[section]

\newtheorem{Theorem}[thm]{Theorem}

\newtheorem{Lemma}[thm]{Lemma}

\newtheorem{Proposition}[thm]{Proposition}
\newtheorem{Corollary}[thm]{Corollary}

\newtheorem{ccote}[thm]{}

{\theorembodyfont{\rmfamily} \theoremstyle{plain}

\newtheorem{Remark}[thm]{Remark}

\newcommand{\preu}{\noindent{\sc Proof: \ }}

\newcommand{\ppreu}[1]{\noindent {\sc Proof of #1: \ }}
\newcommand{\fl}[1]{\buildrel{#1}\over{\longrightarrow}}

\newcommand{\cqfd}{\unskip\kern 6pt\penalty 500
\raise -2pt\hbox{\vrule\vbox to10pt{\hrule width
 8pt\vfill\hrule}\vrule}\smallskip}

\newcommand{\bbr}{{\mathbb{R}}}
\newcommand{\bbc}{{\mathbb{C}}}

\newcommand{\bbz}{{\mathbb{Z}}}

\newcommand{\bbn}{{\mathbb{N}}}

\newcommand{\calb}{{\mathcal B}}
\newcommand{\calc}{{\mathcal C}}

\newcommand{\calh}{{\mathcal H}}

\newcommand{\pcirc}{\kern .7pt {\scriptstyle \circ} \kern 1pt}

\newcommand{\hfl}[1]{\buildrel{#1}\over{\longrightarrow}}
\newcounter{exo}

\newcommand{\mun}{{-1}}
\newcommand{\sk}[1]{\vskip #1 mm}

\newcommand{\onto}{\to\kern-6.5pt\to}
\newcommand{\into}{\hookrightarrow}

\newcommand{\algg}{\mathfrak{g}}

\newcommand{\proref}[1]{Proposition~\ref{#1}}

\newcommand{\lemref}[1]{Lemma~\ref{#1}}
\newcommand{\corref}[1]{Corollary~\ref{#1}}
\newcommand{\thref}[1]{Theorem~\ref{#1}}

\newcommand{\secref}[1]{\S~\ref{#1}}

\newcommand{\csp}{conjugation space}

\newcommand{\hfra}{$H^*$-frame}
\newcommand{\ef}{equivariantly formal}

\newcommand{\R}{\mathbb{R}}

\newcommand{\kappafix}{\kappa_{\scriptscriptstyle\rm fix}}
\newcommand{\sigmafix}{\sigma_{\scriptscriptstyle\rm fix}}
\newcommand{\bkappafix}{\kappa_{\scriptscriptstyle\rm fix}}
\newcommand{\bsigmafix}{\sigma_{\scriptscriptstyle\rm fix}}
\newcommand{\rhofix}{\rho_{\scriptscriptstyle\rm fix}}
\newcommand{\rfix}{r_{\scriptscriptstyle\rm fix}}

\newcommand{\scr}{\scriptscriptstyle}

\newcommand{\dfn}[1]{{\it #1}}
\renewcommand{\:}{\colon}
\newcommand{\lt}[1]{\ell t_{#1}}

\newcommand{\skel}[2]{{\rm Sk}^{#1}_{#2}}
\newcommand{\sq}{{\rm Sq}}

\renewcommand{\emph}{\it}

\title{Conjugation spaces and edges \\ of compatible torus actions}
\author{Jean-Claude HAUSMANN and Tara HOLM}
\date{\ddate}

\begin{document}
\maketitle
\tableofcontents

\section{Introduction}

Duistermaat introduced the {\emph real locus} of a Hamiltonian manifold \cite{Du}.
In this and in others' subsequent works \cite{BGH,Go,GH,HH,Ho,OS,Sd}, it has been
shown that many of the techniques developed in the symplectic category can be
used to study real loci, so long as the coefficient ring is restricted to the
integers modulo $2$.  As we will see, these results seem not necessarily to
depend on the ambient symplectic structure, but rather to be topological in
nature.  This observation prompts the definition of {\emph conjugation space}
in \cite{HHP}.
We now give a brief survey of the results in symplectic geometry that
motivated the definition of a {\emph conjugation space}.

A {\emph symplectic manifold} is a manifold $M$ together with a $2$-form
$\omega\in \Omega^2(M)$ that is closed ($d\omega=0$) and non-degenerate (for
each $X\in T_pM$ there exists $Y\in T_pM$ such that $\omega_p(X,Y)\neq 0$).
Let $G$ be a compact Lie group acting on $M$ preserving $\omega$, $\algg$ the
Lie algebra of $G$, $\algg^*$ its dual, and $\langle
\cdot,\cdot\rangle:\algg^*\times\algg\to \R$ the natural pairing. For each
$X\in\algg$, we let $X^{\#}$ denote the vector field on $M$ generated by the
one-parameter subgroup $\exp(tX)$.  We say that that the $G$-action on $M$ is
{\emph Hamiltonian} if there is a {\emph moment map}
$$
\Phi: M\to\algg^*
$$
that satisfies
\begin{enumerate}
\item $\imath_{X^\#}\omega = d\langle \Phi,X\rangle$ for all $X\in\algg$; and
\item $\Phi$ is equivariant with respect to the given $G$ action on $M$ and the coadjoint action on $\algg^*$.
\end{enumerate}
The function $\Phi_X = \langle \Phi,X\rangle$ is called the {\emph Hamiltonian function} for the vector field $X^\#$.

When $G = T$ is a torus, the second condition on $\Phi$ requires that it be a $T$-invariant map.  In this special case, we have
\begin{thm}[\cite{At:convexity},\cite{GS:convexity}]\label{thm:convex}
If $M$ is a compact Hamiltonian $T$-space, then $\Phi(M)$ is a convex
polytope.  It is the convex hull of $\Phi(M^T)$, the images of the
$T$-fixed points.
\end{thm}
More generally, there are results first of Kirwan and then many others for non-abelian groups.

A by-product of Atiyah's proof of Theorem~\ref{thm:convex} is that any of the
Hamiltonian functions $\Phi_X$ is a perfect Morse function on $M$, in the sense of Bott, and for generic $X$, the critical set is $M^T$.   More precisely,
\begin{eqnarray}\label{1}
H^*(M;\R) & = & \sum_{i=1}^NH^{*-d_i}(F_i;\R),
\end{eqnarray}
where the $F_i$ are the connected components of $M^T$ and $d_i$ is
the Morse-Bott index of $F_i$. This statement is also true over $\bbz$ provided that the cohomology of each $F_i$ is torsion-free, or when the torus action satisfies some additional hypotheses.

Duistermaat introduced the concept of { real locus} to this framework \cite{Du}.  Let $M$ be a Hamiltonian $T$-space, and  $\tau: M\to M$ an anti-symplectic involution that is {\it compatible} with the action; that is, it satisfies
 $$
 \tau(t\cdot p) = t^{-1}\cdot \tau(p),
 $$
 for all $t\in T$ and $p\in M$.  Then if it is non-empty, the submanifold
 $M^\tau$ of $\tau$-fixed points is a Lagrangian submanifold of $M$ called the
 {\it real locus} of the involution.  The primary example of such an
 involution is the one induced by complex conjugation on a complex projective
 variety defined over $\bbr$.  For example, if $M = \bbc P^n$ equipped with
 the Fubini-Study symplectic form and the standard $T^n$ action, then the real
 locus for complex conjugation consists of the real points $\bbr P^n$, whence
 the name real locus.  The main results in \cite{Du} generalize
 Theorem~\ref{thm:convex} and Atiyah's Morse theoretic results.

 \begin{thm}[\cite{Du}]
If $M$ is a compact Hamiltonian $T$-space, and $\tau$ a compatible involution, then
\begin{enumerate}
\item The real locus has full moment image: $\Phi(M^\tau) = \Phi(M)$ is a convex polytope; and
\item Components $\Phi_X$ of the moment map are perfect Morse functions on $M^\tau$, in the sense of Bott, and for generic components the critical set is $M^\tau\cap M^T$, when the coefficients are taken in $\bbz_2$.
\end{enumerate}
\end{thm}

We have the following immediate corollary, a real locus version of
Equation~\eqref{1}, that generalizes classical results on real projective
space and real flag varieties.
\begin{Corollary}\label{2-cor}
If $M$ is a compact Hamiltonian $T$-space, and $\tau$ a compatible involution, then
\begin{eqnarray}\label{2}
H^*(M^\tau;\bbz_2) & = & \sum_{i=1}^NH^{*-\frac{d_i}{2}}((F_i)^\tau;\bbz_2) \, .
\end{eqnarray}
where the $F_i$ are the connected components of $M^T$ and $d_i$ is the Morse-Bott index of $F_i$ (in $M$).
\end{Corollary}

Duistermaat's work began a flurry of activity on properties of real
loci. We provide a brief account here; a more detailed record is
available in \cite{Sj}.
Davis and Janusz\-kiewicz studied the real loci of toric varieties in
their own right \cite{DJ}, independent of Duistermaat's work.
The first author and Knutson analyze a large class of
examples of real loci in their account of planar and spacial polygon spaces
\cite{HK}. O'Shea and Sjamaar generalized Kirwan's non-abelian convexity
results to real flag manifolds and real loci \cite{OS}.  This has recently
been extended by Goldberg \cite{Go}.

Schmid and independently Biss, Guillemin and the second author generalized
\eqref{2} to the equivariant setting: the idempotents $T_2 = \{ t\in T\ |\ t^2
= 1\}$ act on the real locus, and many results in $T$-equivariant symplectic
geometry may be generalized to $T_2$ equivariant geometry of real loci (with
coefficients restricted to $\bbz_2$) \cite{BGH,Sd}.
This work yields an explicit description of the $T_2$-equivariant cohomology
for the fixed set of the Chevalley involution on certain coadjoint orbits, and
on the real locus of a toric variety, using localization methods.  These
results were strengthened to include the fixed set of the Chevalley involution
on all coadjoint orbits in \cite{HHP}.

Following this, Goldin and the second author \cite{GH} proved that there is a
natural involution on an abelian symplectic reduction of a symplectic manifold
with involution.  Moreover, the $T_2$ equivariant cohomology of the original
real locus surjects onto the ordinary cohomology of the real locus of the
symplectic reduction.  This includes a comprehensive description of toric
varieties and their real loci from yet a third perspective.

In all of these papers, a common theme is that there is a degree-halving isomorphism
$$
H^{2*}(M;\bbz_2) \to H^*(M^\tau;\bbz_2).
$$
As we now describe, this can be seen as part of a purely topological
framework, that of a {\emph conjugation space}, introduced in \cite{HHP}.
The remainder of the article is organized as follows.
In Section~\ref{sec:conj}, we review the definitions and properties of \csp s.  Our main theorem
gives a criterion for recognizing when a topological space is a conjugation space; this
is stated in Section~\ref{S.CS&1}, along with two noteworthy corollaries. We then prove some basic facts in
Section~\ref{S.preli}, and prove the main theorems in Section~\ref{sec:mainproofs}.

\medskip

\noindent {\bf NOTE:} For the remainder of the paper, the cohomology is taken
coefficients in the field  $\bbz_2$: $H^*(X)=H^*(X;\bbz_2)$.

\section{A review of \csp s}\label{sec:conj}

Let $X$ be a $G$-space $X$ for a topological group $G$.
The equivariant cohomology $H^*_G(X)$ is defined as the (singular) cohomology of
the Borel construction:
$$
H^*_G(X)=H^*(X\times_G BG) \, .
$$
Hence, $H^*_G(X)$ is a $H^*(BG)$-algebra.
When $G=C$ is the group of order two, $BC=\bbr P^\infty$ and
$H^*(BC)=\bbz_2[u]$, with $u$ in degree $1$. Thus,
$H^*_C(X)$ is a $\bbz_2[u]$-algebra.

Let $\tau$ be a continuous involution on a space $X$.
Let $\rho\: H^{2*}_C(X)\to H^{2*}(X)$ and
$r\: H^{*}_C(X)\to H^{*}_C(X^\tau)$ be the restriction homomorphisms,
where $C=\{{\rm id},\tau\}$.

A {\it cohomology frame} or {\it \hfra}\
for $(X,Y)$ is a pair $(\kappa,\sigma)$, where
\renewcommand{\labelenumi}{(\alph{enumi})}
\begin{enumerate}
\item $\kappa\:  H^{2*}(X)\to H^{*}(X^\tau)$ is
an additive isomorphism dividing the degrees in half; and
\item $\sigma\:  H^{2*}(X)\to H^{2*}_C(X)$ is
an additive section of $\rho$.
\end{enumerate}
Moreover, $\kappa$ and $\sigma$ must satisfy the {\it conjugation equation}
\begin{equation}\label{defalceq}
r\pcirc\sigma(a) = \kappa(a)u^m + \lt{m}
\end{equation}
for all $a\in H^{2m}(X)$ and all $m\in\bbn$, where $\lt{m}$
denotes any polynomial in the variable $u$ of degree less than $m$.
An involution admitting a \hfra\ is called a \dfn{conjugation}.
An even cohomology space (i.e. $H^{odd}(X)=0$) together with a conjugation is called a \dfn{\csp}.
Conjugation spaces were introduced in \cite{HHP} and studied further in \cite{FP}
and \cite{Ol}. The main examples of conjugations are given the complex conjugation in flag manifolds,
the Chevalley involution in coadjoint orbit of compact Lie groups and other natural involutions, e.g. on toric manifolds
or polygon spaces.
Here below are some important properties of \csp s.

\begin{enumerate}
\item If $(\kappa,\sigma)$ is \hfra, then $\kappa$ and $\sigma$ are ring homomorphisms \cite[Theorem~3.3]{HHP}.
The ring homomorphism $\kappa$ also commutes with the Steenrod squares:
$\kappa\pcirc\sq^{2i}=\sq^i\pcirc\kappa$, \cite[Theorem~1.3]{FP}.
\item \hfra s are natural for $\tau$-equivariant maps \cite[Prop.3.11]{HHP}. In particular, if an involution
admits an \hfra, it is unique \cite[Cor.3.12]{HHP}.
\item For a conjugate-equivariant complex vector bundle $\eta$ (``real bundle''
in the sense of Atiyah) over a \csp\ $X$, the isomorphism $\kappa$
sends the total Chern class of $\eta$ onto the total Stiefel-Whitney class of
its fixed bundle.
\end{enumerate}

Duistermaat's \corref{2-cor} admits the following generalization, proved in
\cite[Theorem~8.3]{HHP}.

\begin{Theorem}\label{HamilCS}
Let $M$ be a compact symplectic manifold equipped with a
Hamiltonian action of a torus $T$ and with
a compatible smooth anti-symplectic involution
$\tau$. If $M^T$ is a \csp, then $M$ is a \csp.
\end{Theorem}

The proof of \thref{HamilCS} involves properties of conjugations compatible with $T$-actions
which are interesting by their own.
The involution $g\mapsto g^\mun$ on the torus $T$ induces an involution on $ET$.
Using this involution together with $\tau$, we get an involution on
$X\times ET$ which descends to an involution, still called $\tau$, on $X_T$.
To a torus $T$ is associated its $2$-torus, i.e. the set of idempotnt elements of $T$:
$$
T_2=\{g\in T\mid g^2=1\} \, .
$$
The compatibility implies that $T_2$ acts on $X^\tau$.
The following lemma is proved in \cite[Lemma~7.3]{HHP}

\begin{Lemma}\label{CBTfix}
$(X_T)^\tau=(X^\tau)_{T_2}$. \cqfd
\end{Lemma}

The following theorem is proved in \cite[Theorem~7.5]{HHP}.

\begin{Theorem}\label{strcomp}
Let $X$ be a \csp\ together with a compatible action of
a torus $T$. Then the involution induced on $X_T$ is a conjugation.
\cqfd
\end{Theorem}

Using Lemma~\ref{CBTfix}, one gets the following corollary of
Theorem~\ref{strcomp}.

\begin{Corollary}\label{srtcomp-cor}
Let $X$ be a \csp\ together with an involution and a compatible
$T$-action. Then there is a ring isomorphism
$$\bar\kappa\: H^{2*}_T(X)\hfl{\approx}{} H^{*}_{T_2}(X^\tau). \cqfd$$
\end{Corollary}

\section{Conjugation spaces and $1$-skeleta}\label{S.CS&1}
We now state our new results. They consist of criteria to
determine that an involution $\tau$ is a conjugation, in the case
where $\tau$ is compatible with an action of a torus $T$.
The conditions are on the equivariant $1$-skeleton of the
action of $T$ on $X$ and of the inherited action of the associated $2$-torus $T_2$.

Let $X$ be a topological space, together with a continuous action
of a group $G$, where $G$ is a torus or a $2$-torus
(finite elementary abelian $2$-group).
We define the {\it $G$-equivariant $i$-skeleton
$\skel{G}{i}(X)$
of the $G$-action on $X$} to be
\begin{equation}\label{defeqsk}
\skel{G}{i}(X)=\{x\in X\mid {\rm codim\,}(G_x\subset G)\leq i\},
\end{equation}
where $G_x$ denotes the $G$-isotropy group of $x$.
In~\eqref{defeqsk}, the ``codimension'' is interpreted as
the codimension of a manifold if $G$ is a torus, and
the codimension of a $\bbz_2$-vector subspace if $G$ is
a $2$-torus (and hence isomorphic to a $\bbz_2$-vector space).
In particular,
$\skel{G}{0}(X)$ is equal to the subspace $X^G$ of fixed points.
An {\it edge} (of the $G$-action) is the closure of a connected component
of the set $\skel{G}{1}(X)-\skel{G}{0}(X)$.

Let $T$ be a torus and $T_2$ the subgroup of idempotents. A $T$-action on a space $X$ induces a $T_2$-action on $X$ that satisfies $\skel{T}{i}(X)\subset \skel{T_2}{i}(X)$. For example, $X^T\subset X^{T_2}$.

A continuous action of a topological group $G$ on a space $X$ is
called \dfn{good} if $X$ has the $G$-equivariant homotopy type
of a finite $G$-CW-complex. For instance, a smooth action of
a compact Lie group on a closed manifold is good.
A continuous involution $\tau$ is
called \dfn{good} if the corresponding action of the cyclic group
$C=\{{\rm id},\tau\}$ is good.

Let $X$ be a topological space, and let $\tau$ be continuous involution on $X$
that is compatible with a continuous action of a torus $T$.
Then the involution $\tau$
preserves the $T$-equivariant skeleta and sends an each edge to a (possibly different) edge.
Moreover, the real locus $X^\tau=X^C$ inherits an action of $T_2$.
Our main results are the following.

\begin{Theorem}[Main Theorem]\label{Prop1sk}
Let $X$ be an even cohomology space, together with a good involution $\tau$
which is compatible with a good action of a torus $T$.
Suppose that
\begin{enumerate}
\renewcommand{\labelenumi}{(\alph{enumi})}
\item $(X^T,\tau)$ is a conjugation space.
\item each edge of the $T$-action is preserved by $\tau$
and  is a conjugation space.
\item  $\skel{T}{i}(X)=\skel{T_2}{i}(X)$ for $i=0,1$.
\end{enumerate}
Then $X$ is a conjugation space.
\end{Theorem}

Recall that a $T$-action on a space $X$ is called a \dfn{GKM}
action if each edge is a $2$-sphere upon which $T$ acts by rotation around
some axis, via a non-trivial character $T\to S^1$. One consequence of this assumption is that
$X^T$ is discrete.

\begin{Corollary}\label{C.GKM}
Let $X$ be an even cohomology space, together with a good involution $\tau$
which is compatible with a good GKM action of a torus $T$, satisfying
$\skel{T}{i}(X)=\skel{T_2}{i}(X)$ for $i=0,1$.
Suppose that $\tau$ acts trivially on $X^T$ and preserves each edge.
Then $X$ is a conjugation space.
\end{Corollary}

\begin{Corollary}\label{C.edgeHam}
Let $X$ be an even cohomology space, together with a good involution $\tau$
which is compatible with a good action of a torus $T$, satisfying
$\skel{T}{i}(X)=\skel{T_2}{i}(X)$ for $i=0,1$.
Suppose that
\begin{enumerate}
\renewcommand{\labelenumi}{(\alph{enumi})}
\item $(X^T,\tau)$ is a conjugation space.
\item each edge of the $T$-action is preserved by $\tau$ and is
a Hamiltonian $T$-manifold on which $\tau$ acts smoothly and is
anti-symplectic.
\end{enumerate}
Then $X$ is a conjugation space.
\end{Corollary}

See \secref{S.comm} for comments about the condition $\skel{T}{i}(X)=\skel{T_2}{i}(X)$ for $i=0,1$.

\section{Preliminaries}\label{S.preli}

This section is devoted to the proof of \thref{Prop1sk} and of Corollaries~\ref{C.GKM}
and~\ref{C.edgeHam}. We begin with some preliminaries.

\begin{ccote}\label{compatibility} Compatibility. \  \rm
Let $X$ be a topological space endowed with a continuous involution $\tau$
which is compatible with a continuous action of a torus $T$.
The involution $\tau$ then induces an involution on on the fixed point set $X^T$.
In addition, the associated $2$-torus $T_2$ of $T$ acts on $X^\tau$
and $X^\tau\cap X^T\subset (X^\tau)^{T_2}$.
Condition~(c) of \thref{Prop1sk} will play an important role.

\begin{Lemma}\label{L.XTXT2}
Suppose that $\skel{T}{i}(X) = \skel{T_2}{i}(X)$. Then
$\skel{T}{i}(X)^\tau=\skel{T_2}{i}(X^\tau)$.
\end{Lemma}

\preu One has
\begin{equation}\label{T=T21stcon}
\skel{T}{i}(X)^\tau = \skel{T}{i}(X)\cap X^\tau\subset
\skel{T_2}{i}(X^\tau) = X^\tau\cap \skel{T_2}{i}(X) =
X^\tau\cap\skel{T}{i}(X) \ ,
\end{equation}
which implies that $\skel{T}{i}(X)^\tau=\skel{T_2}{i}(X^\tau)$.
\cqfd

\end{ccote}

\begin{ccote}\label{eqform} Equivariantly formal spaces.\ \rm
Let $X$ be a space with a continuous action
of a compact Lie group $G$.
First introduced in \cite{GKM} for $G$ a torus and complex coefficients, the notion of
\textit{equivariant formality} was developed for other coefficients where it
is more subtle, see \cite[(2.3)]{HHP} and \cite[\S~8]{Fr}.
A $G$-space $X$ is {\it \ef\ (over $\bbz_2$)} if
the map $X\to EG\times_G X$ is totally non-homologous to zero,
that is the restriction homomorphism
$j^*\:H^*_G(X)\to H^*(X)$ is surjective.
A space $X$ with an involution $\tau$ is called $\tau$-\ef\ if
it is $C$-\ef\ for $C=\{{\rm id},\tau\}$.
The following results are classical
but may be not found in the literature with exactly
our hypotheses. Let $R=H^*_G(pt)$; the map $X_G\to BG$ gives
a ring homomorphism $p^*\:R\to H^*_G(X)$, making $H^*_G(X)$
an $R$-module.

\begin{Proposition}\label{Peqform}
The following conditions are equivalent:
\begin{enumerate}
\renewcommand{\labelenumi}{(\roman{enumi})}
\item\label{eqfo}
$X$ is an equivariantly formal $G$-space.
\item\label{d20}
The group $G$ acts trivially on
$H^*(X)$ and the Serre spectral sequence
for the cohomology of the fibration
$X\to EG\times_G X\to BG$ collapses at the term $E_2$.
\item\label{free}
The group $G$ acts trivially on
$H^*(X)$ and $H^*_G(X)$ is a free $R$-module.
\item there is an additive homomorphism
$\sigma\:H^*(X)\to H^*_G(X)$ such that $j^*\pcirc\sigma={\rm id}$
and $p^*\otimes\sigma\:R\otimes H^*(X)\to H^*_G(X)$
is an isomorphism of $R$-modules.
\item \label{tensor}
The ring homomorphism $H^*_G(X)\to H^*(X)$
descends to a ring isomorphism
$H^*_G(X)\otimes_{R}\bbz_2\hfl{\approx}{} H^*(X)$.
\end{enumerate}
\end{Proposition}

\preu
This proof is for ${\rm mod}\,2$-cohomology, but it works for
the cohomology with coefficients in any field.
\sk{1}\noindent\it (i) is equivalent to (ii):\rm\
the ring homomorphism
$j^*:H^*_G(X)\to H^*(X)$ is the composition:
\begin{eqnarray}\label{eqfo1}
H^*_G(X)\onto E_\infty^{0,*}\subset E_2^{0,*} & = &
H^0(BG;H^*(X))\\
& = & H^*(X)^{G}\subset H^*(X). \nonumber
\end{eqnarray}
If these inclusions are equalities, then
$j^*$ is onto, which shows that (ii)
implies (i). Conversely, if $j^*$ is onto,
this shows that
$H^*(X)^{G}=H^*(X)$ and
$E_\infty^{0,*}=E_2^{0,*}$.
As the differentials are morphisms of $R$-modules,
this implies that
$E_\infty^{*,*}=E_2^{*,*}$ (see \cite[p. 148]{McC}).
Hence (i) implies (ii).

\sk{1}\noindent\it (i) implies (iii) and (iv) :\rm\
As $j^*$ is surjective, there exists a $\bbz_2$-linear
section $\sigma$ of $j^*$. We already showed that (i) implies
that the $G$-action on $H^*(X)$ is trivial.
As $G$ is a compact Lie group,
$H^p_G(pt)$ is a finite dimensional
$\bbz_2$-vector spaces for all $p$.
The Leray-Hirsch theorem \cite[Thm 5.10]{McC} then implies
that $H^*_G(X)$ is a free $R$-module
with basis $\sigma(\calb)$, where $\calb$ is a $\bbz_2$-basis of
$H^*(X)$. This implies (iii) and (iv).

\sk{1}\noindent\it (iii) implies (i):\rm\
as $G$ acts trivially on $H^*(X)$,
the term $E_2^{*,*}$ is isomorphic to $R\otimes H^*(X)$
as a bigraded $R$-module. This implies that the
kernel of $j^*$ is equal to $I\cdot H^*_G(X)$,
where $I$ is the ideal of $R$ of elements of positive degree.
Suppose that $H^*_G(X)$ is the free $R$-module with
some basis $\calc$. As $R/I=\bbz_2\otimes_R R\approx\bbz_2$,
the image of $j^*$ can be identified with $\bbz_2$-vector space with
basis $\calc$. Denote by $\calc_{s}$ the subset of $\calc$
of elements of degree $\leq s$.

Suppose, by induction on $q$, that
$j_*:H^q_G(X)\to H^q(X)$ is surjective for $q\leq k$
(true for $k=0$).
If there is $a\in H^k(X)$ which is not
in the image of $j^*$, then $d_r(a)\neq 0$ for some
differential $d_r:H^k(X)\to R\otimes H^{k-r+1}(X)$.
Therefore, there are elements $a_1,\dots,a_m\in \calc_{k-r-1}$,
and $r_1,\dots,r_m\in R$ with $\sum r_ia_i=d_r(a)$
in $E_{\infty}^{*,*}$. This means that $\sum r_ia_i\in H^*_G(X)$
is a $R$-linear combination
of elements of $\calc_{k-r-2}$. Such a relation would contradict
the fact that $\calc$ is a basis of $H^*_G(X)$.

\sk{1}\noindent\it (iv) implies (v):\rm\
the homomorphism $j^*\pcirc p^*: R\to H^*(X)$ coincides
with the projection $R\to\bbz_2\otimes_R R = \bbz_2$. Therefore,
$j^*$ factors through a ring homomorphism
$\bar j^* : \bbz_2\otimes_R H^*_G(X)\to H^*(X)$. On the other hand,
$j^*\pcirc\sigma={\rm id}$. Hence, one has a commutative
diagram
$$\begin{array}{cccccccc}
R\otimes H^*(X) & \hfl{}{} & \bbz_2\otimes_R(R\otimes H^*(X))
& = &  H^*(X) \\
p^*{\scriptscriptstyle\otimes}\,\sigma\downarrow\,\approx && \downarrow\approx
&& \downarrow\,\, {\rm id} \\
H^*_G(X) & \hfl{}{} &  \bbz_2\otimes_R H^*_G(X)
& \hfl{\bar j^*}{} &  H^*(X) \ ,
\end{array}$$
which proves that $\bar j^*$ is an isomorphism.

\sk{1}\noindent\it (v) implies (i):\rm\  this implication is trivial.\cqfd

\begin{Proposition}\label{Peqform2}
Let $X$ be a good $G$-space which
is equivariantly formal over $\bbz_2$.
Suppose that one of the following hypotheses holds:
\begin{enumerate}
\renewcommand{\labelenumi}{(\alph{enumi})}
\item $G$ is a torus and $X^G=X^{G_2}$.
\item $G$ is a $2$-torus.
\end{enumerate}
Then the restriction homomorphism
$H^*_G(X)\to H^*_G(X^G)$ is injective.
\end{Proposition}

\begin{Remark}\rm In Case (b), Proposition
\ref{Peqform2} is false without the assumption $X^G=X^{G_2}$.
For example, consider the $G=S^1$ action on $X=S^2\subset\bbc\times\bbr$
by $g(z,t)=(g^2z,t)$. This has $X^G=\{(0,\pm 1)\}$. Let
$U_+=X-\{(0,-1)\}$ and $U_-=X-\{(0,-1)\}$. The intersection
$U_+\cap U_-$ is $G$-homotopy equivalent to the homogeneous
space $G/G_2$ and then $H^*_G(U_+\cap U_-)=H^*(BG_2)$. The
Mayer-Vietoris sequence for $(X,U_+,U_-)$ then gives
$$
0 \to  H^1(BG_2) \to  H^2_G(X) \to H^2_G(X^G)
$$
and $H^1(BG_2)=\bbz_2$.
\end{Remark}

\ppreu{Proposition \ref{Peqform2}}
Let $R_{(0)}$ be the field of fractions of $R$,
that is $R$ localized at $S=R-\{0\}$.
By our assumptions, the multiplicative set $S$ is central in $R$.
Let
$$
X^S=\{x\in X\mid H^*(BG)\to H^*(BG_x) \hbox{ is injective} \} \, ,
$$
where $G_x$ is the isotropy group of $x$.
The localization theorem (\cite[Thm.3.1.6]{AP}, \cite[Thm.3.7]{Al})
asserts that the inclusion $X^S\subset X$ induces an isomorphism of
$R_{(0)}$-vector spaces
\begin{equation}\label{locaa}
S^\mun H^*_G(X)\approx S^\mun H^*_G(X^S).
\end{equation}
In Case (b), if $G_x$ is a proper subgroup of $G$, then
$H^2(BG)\to H^2(BG_x)$ is not injective; hence $X^S = X^G$.
For Case (a), we use that, for each $x\in X$, there is an isomorphism
$\psi_x\:G\hfl{\approx}{} (S^1)^m$ such that
$\psi_x(G_x)=C_1\times\cdots\times C_m$, where $C_j$ is a subgroup
of $S^1$. In order to have
$H^2(BG)\to H^2(BG_x)$ injective, each $C_j$ should
be either $S^1$ or a finite cyclic group of even order.
Then $X^G\subseteq X^S \subseteq X^{G_2}=X^G$.
Hence, in all cases, we have proved that
\begin{equation}\label{loca}
S^\mun H^*_G(X)\to S^\mun H^*_G(X^G)
\end{equation}
is an isomorphism.
Therefore, $\ker(H^*_G(X)\to H^*_G(X^G))$ is the
$R$-torsion of $H^*_G(X)$.  But the $R$-torsion vanishes because $H^*_G(X)$
is a free $R$-module by Proposition~\ref{Peqform}.\cqfd

\begin{Proposition}\label{Peqform3}
Let $X$ be a good $G$-space.
Suppose that Conditions (a) or (b) of Proposition~\ref{Peqform2}
are satisfied.
Then $X$ is an \ef\ $G$-space over $\bbz_2$
if and only if
${\rm dim}_{\bbz_2} H^*(X) = {\rm dim}_{\bbz_2} H^*(X^G)$.
\end{Proposition}

\preu
As in the proof of Proposition \ref{Peqform2},
consider $S=R-\{0\}$ and $R_{(0)}=S^\mun R$.
We apply $S^\mun$ to the
terms of the Serre spectral sequence, following
\cite[proof of Cor.3.10]{Al}.
When $G$ is a torus, it
acts trivially on $H^*(X)$, which implies that
$E_2^{*,*}\approx H^*(BG;H^*(X))$ as $R$-module and
there is an isomorphism of  $R_{(0)}$-vector spaces
$R_{(0)}\otimes_{\bbz_2} H^*(X)\approx S^\mun E_2$.
Therefore, using equation~\eqref{loca}, we get
\begin{eqnarray}\label{eqform3}
\dim_{\bbz_2} H^*(X) & = &  {\rm dim\,}_{R_{(0)}}(S^\mun E_2)
\nonumber \\
 & \geq & \dim_{R_{(0)}}(S^\mun E_\infty) \nonumber \\
 & = & \dim_{R_{(0)}}(S^\mun H^*_G(X)) \\
 & = & \dim_{R_{(0)}}(S^\mun H^*_G(X^G))\nonumber \\
 & = & \dim_{\bbz_2}(H^*(X^G)).\nonumber
\end{eqnarray}
By Proposition \ref{Peqform},
the inequality in equation~\eqref{eqform3} is an equality
if and only if $X$ is \ef.
Finally, when $G$ is a $2$-torus,
Proposition~\ref{Peqform3} follows from
\cite[Thm\,3.10.4]{AP}.\cqfd
\end{ccote}

\begin{ccote}\rm
We shall need the following two lemmas, first proved by
Chang and Skjelbred for rational cohomology and torus action~\cite{CS}.

\begin{Lemma}\label{L.CS2}
Let $X$ be a space endowed with a good action of a $2$-torus $G$.
Suppose that $X$ is $G$-\ef.
Then the  restriction homomorphisms on the $\rm mod\,\, 2$-cohomology
$H^*_G(X)\to H^*_G(X^G)$
and $H^*_G(\skel{G}{1}(X))\to H^*_G(X^G)$ have same image.
\end{Lemma}

\preu
Using the equivalence (i) $\Leftrightarrow$ (iii) in
\lemref{Peqform}, we know that $H^*_G(X)$ is a free $H^*_G(pt)$-module.
By \cite[Corollary p.~63]{Hs}, the homomorphism
$$
H^*_G(X,X^G) \to H^*_G(\skel{G}{1}(X),X^G)
$$
is injective. The $H^*_G$-sequences of the pairs
$(X,X^G)$ and $(\skel{G}{1}(X),X^G)$ are part of a commutative diagram
$$
\begin{array}{c}{\xymatrix@C-3pt@M+2pt@R-4pt{%
0
\ar[r]  &
H^*_G(X) \ar[d]
\ar[r]^(0.50){}  &
H^*_G(X^G) \ar[d]^(0.50){=} \ar[r]
&
H^{*+1}_G(X,X^G) \ar@{>->}[d] \ar[r] & 0
\\
&
H^*_G(\skel{G}{1}(X)) \ar[r]  &
H^*_G(X^G) \ar[r]  &
H^{*+1}_G(\skel{G}{1}(X),X^G) \ar[r] & 0 \ .
}}\end{array}
$$
Therefore, the injectivity of the last vertical arrow implies the lemma.
\cqfd

The following lemma follows from
\cite[Theorem~2.1]{FP3}.

\begin{Lemma}\label{L.CS}
Let $X$ be a space endowed with a good action of a torus $T$.
Suppose that $X$ is $T$-\ef\ and that
$\skel{T}{i}(X) = \skel{T_2}{i}(X)$ for $i=0,1$.
Then the  restriction homomorphisms on the $\rm mod 2$-cohomology
$H^{*}_T(X)\to H^{*}_T(X^T)$
and $H^{*}_T(\skel{T}{1}(X))\to H^{*}_T(X^T)$ have same image.
\cqfd
\end{Lemma}
\end{ccote}


\section{Proof of the main results}\label{sec:mainproofs}
We begin with the proof of \thref{Prop1sk}.  Recall that we are working with cohomology with $\bbz_2$ coefficients.  In what follows, $\dim$ denotes $\dim_{\bbz_2}$.

\begin{ccote}\label{Xef}
$X$ is is $\tau$-equivariantly formal over $\bbz_2$ and
$X^\tau$ is is $T_2$-equivariantly formal over $\bbz_2$. \rm
Being an even-cohomology space, $X$ is $T$-\ef.
By Hypothesis (a) of \thref{Prop1sk} and by \lemref{L.XTXT2},
we have
\begin{eqnarray*}
\dim H^{*}(X) &=& \dim H^{*}(X^T) = \dim H^{*}((X^T)^\tau)
\\ &=& \dim H^{*}((X^\tau)^{T_2})
\leq \dim H^{*}(X^\tau) \leq \dim H^{*}(X) \, ,
\end{eqnarray*}
which implies that
\begin{equation}
\dim H^{*}(X^\tau) =  \dim H^{*}(X)  \ \hbox{ and }
\dim H^{*}((X^\tau)^{T_2}) =\dim H^{*}(X^\tau) \, .
\end{equation}
\end{ccote}

\begin{ccote}\label{XTef} $X_T$ is  $\tau$-\ef.\ \rm
For $G$ a topological group and $k\in\bbn$,
we consider the $G$-principal bundle $G\to E_kG\to B_kG$
obtained as $k$-th step in the Milnor construction.
If $X$ is a $G$-space, the associated bundle with fibre $X$
gives a bundle $X\to X_{G,k}\to B_kG$,
where $X_{G,k}= E_kG\times_G X$.

For a torus $T$ of dimension $n$, $B_kT\approx  (\bbc P^k)^n$.
The involution $\tau(g)=g^\mun$ on $T$ gives an involution
$\tau$ on $B_kT$ which makes $B_kT$ a conjugation space with
$(B_kT)^\tau\approx (\bbr P^k)^n \approx B_kT_2$.

We first prove that $X_{T,k}$ is $\tau$-\ef. As $X$ and $B_kT$
are even cohomology spaces, the spectral sequence of $X\to X_{T,k}\to B_kT$
degenerates at the $E^2$-term and
$H^*(X_{T,k})\approx H^*(X)\otimes H^*(B_kT)$. As a consequence,
$\dim H^*(X_{T,k}) = \dim H^*(X)\cdot \dim H^*(B_kT)<\infty$.
Since $X$ is $\tau$-\ef\ by \ref{Xef}, one has
\begin{equation}\label{XTef-eq1}
\dim H^*(X_{T,k}) = \dim H^*(X) \cdot \dim H^*(B_kT) =
\dim H^*(X^\tau) \cdot \dim H^*(B_k{T_2})
\end{equation}
As $X^\tau$ is $T_2$-\ef\ by \ref{Xef}, the following commutative
diagram
$$
\begin{array}{c}{\xymatrix@C-3pt@M+2pt@R-4pt{%
H^*((X^\tau)_{T_2}) \ar@{>>}[r]^(0.50){\rho^\tau_{T_2}} \ar[d]^(0.50){}  &
H^*(X^\tau) \ar[d]^(0.50){=}  \\
H^*((X^\tau)_{T_2,k}) \ar[r]^(0.50){\rho^\tau_{T_2,k}} &
H^*(X^\tau)
}}\end{array}
$$
shows that $\rho^\tau_{T_2,k}$ is surjective and thus
$H^*((X^\tau)_{T_2,k})\approx H^*(X^\tau)\otimes H^*(B_kT_2)$. As
in \lemref{CBTfix}, one has $(X_{T,k})^\tau=(X^\tau)_{T_2,k}$, thus
\begin{equation}\label{XTef-eq2}
\dim H^*((X_{T,k})^\tau) = \dim H^*((X^\tau)_{T_2,k})=
\dim H^*(X^\tau) \cdot \dim H^*(B_kT_2) \, .
\end{equation}
Putting \eqref{XTef-eq1} and \eqref{XTef-eq2} together gives
$\dim H^*(X_{T,k})=\dim H^*((X_{T,k})^\tau)$, and
with \proref{Peqform3}, this implies that $X_{T,k}$ is \ef.

Now given $n\in\bbn$, there exists $k\in\bbn$ such that
$H^n(X_T)\approx H^n(X_{T,k})$. The following commutative diagram
$$
\begin{array}{c}{\xymatrix@C-3pt@M+2pt@R-4pt{%
H^n(X_{T}) \ar[r]^(0.50){\rho} \ar[d]^(0.50){\approx}  &
H^n(X) \ar[d]^(0.50){=}  \\
H^n(X_{T,k}) \ar@{>>}[r]^(0.50){\rho_k} &
H^n(X)
}}\end{array}
$$
shows that $\rho$ is surjective in degree $n$. This can be done for
each $n$, so $X_T$ is \ef.
\end{ccote}

\begin{ccote}\label{constrkt}
Construction of the ring isomorphism $\kappa_{\scr T}\: H^{2*}(X_T)\to H^{*}((X_T)^\tau)$. \rm
By \lemref{CBTfix}, it is equivalent to construct a
ring isomorphism $$\kappa_{\scr T}\: H^{2*}_T(X)\to H^{*}_{T_2}(X^\tau).$$
By \corref{srtcomp-cor}, such an isomorphism
$\bkappafix\:H^{2*}_T(X^T)\to H^{*}_{T_2}((X^T)^\tau)$ exists, since $X^T$ is a
\csp.
As $(X^T)^\tau=(X^\tau)^{T_2}$ by \lemref{L.XTXT2}, we may view $\bkappafix$ as a map from
$H^{2*}_T(X_T)$ to $H^{*}_{T_2}((X^\tau)^{T_2})$. Consider the
following diagram.
\begin{equation}\label{constrdiag}
\begin{array}{c}
\xymatrix{
H^{2*}_T(X) \ar[r]^{q} & H^{2*}_T(X^T)\ar[d]^{\bkappafix}_\approx \\
H^{*}_{T_2}(X^\tau) \ar[r]^(0.45){q^\tau} & H^{*}_{T_2}((X^\tau)^{T_2}) \, .
}
\end{array}
\end{equation}
By Proposition \ref{Peqform2}, the restriction homomorphisms $q$ and $q^\tau$
are injective. Therefore, in order to construct $\kappa_{\scr T}\: H^{2*}_T(X)\to H^{*}_{T_2}(X^\tau)$,
it is enough to show that $A^\tau=\bkappafix(A)$, where
$A={\rm image}(q)$ and
$A^\tau={\rm image}(q^tau)$.  The proof is a diagram chase.

Let $N$ be the $1$-skeleton of $X$ and let $N^\tau=N\cap X^\tau$.
Let $\tilde N$ be the {\it disjoint union} of all the edges of $X$.
There is an obvious quotient map $\tilde N\to N$.
Let $N_0\subset N\cap X^T$ be the points of $N$ having more
than 1 preimage in $\tilde N$, and let $\tilde N_0$ be the points
of $\tilde N$ above $N_0$. Thus $N_0$ is a union of components of $X^T$
and $\tilde N_0\to N_0$ is a disjoint union of trivial coverings.
The various inclusions give a morphism of push-out diagrams
\begin{equation}
\begin{array}{ccccc}
\tilde N_0 &\fl{}{}&\tilde N^T\\
\downarrow  &&  \downarrow \\
N_0  &\fl{}{}& X^T
\end{array}\qquad \longrightarrow \qquad
\begin{array}{ccccc}
\tilde N_0 &\fl{}{}&\tilde N \\
\downarrow  &&  \downarrow \\
N_0  &\fl{}{}& N
\end{array}
\end{equation}
By Hypothesis (a) $X^T$ is a \csp. By \cite[Remark~3.1]{HHP}, $\tau$ preserves each
arc-connected component of $X^T$. Therefore, $N_0$ and $\tilde N_0$ are \csp s.
In the same way, using Hypothesis (b), $\tilde N^T$ and $\tilde N$
are \csp s. The induced morphism on Mayer-Vietoris sequences,
together with the isomorphisms $\kappa$'s
and the fact that that $(X^T)^\tau=(X^\tau)^{T_2}$ (by \lemref{L.XTXT2})
gives	 a three dimensional commutative diagram:
\begin{equation}\label{BD}
{\begin{tiny}
\xymatrix@C-35pt@M+4pt{
    && H^{2*}_T(X^T) \ar'[d][dd]^(0.3){\scriptscriptstyle\approx}_(0.3){\bkappafix}
\ar@{>->}[rr]
    && H^{2*}_T(N_0)\oplus H^{2*}_T(\tilde N^T)
  \ar'[d][dd]^(0.3){\scriptscriptstyle\approx}\ar@{->>}[rr]
   &&  H^{2*}_T(\tilde N_0) \ar[dd]^{\scriptscriptstyle\approx}
\\ 
 & H^{2*}_T(N) \ar[rr]\ar[ur]
 && H^{2*}_T(N_0)\oplus H^{2*}_T(\tilde N)
 \ar[dd]^(0.3){\scriptscriptstyle\approx}\ar[rr]\ar[ur] &&
  H^{2*}_T(\tilde N_0)
  \ar[dd]^(0.3){\scriptscriptstyle\approx}\ar[ur]_{\scriptscriptstyle =}
\\  
H^{2*}_T(X) \ar@/^7mm/[urur]^{q}
\ar[ur] && H^*_{T_2}((X^\tau)^{T_2})
\ar@{>->}'[r][rr] &&
  H^*_{T_2}(N^\tau_0)\oplus H^*_{T_2}((\tilde N^\tau)^{T_2}) \ar@{->>}'[r][rr] &&
  H^*_{T_2}(\tilde N^\tau_0)
\\   
 & H^*_{T_2}(N^\tau) \ar[rr]\ar[ur] && H^*_{T_2}(N^\tau_0)\oplus H^*_{T_2}(\tilde N^\tau)
 \ar[rr]\ar[ur] && H^*_{T_2}(\tilde N^\tau_0) \ar[ur]_{\scriptscriptstyle =}
\\   
H^*_{T_2}(X^\tau) \ar[ur]\ar@/^7mm/[urur]^{q^\tau}
}\end{tiny}}
\end{equation}

The vertical squares commute because of the naturality of the \hfra s
of conjugation spaces.
The quotient maps $\tilde N^T\to X^T$ and $(\tilde N^T)^\tau\to(X^T)^\tau$ admit
continuous sections, so the homomorphisms
$H^{2*}_T(X^T)\to H^{2*}_T(N_0)\oplus H^{2*}_T(\tilde N^T)$ and
$H^*_{T_2}((X^\tau)^{T_2})\to H^*_{T_2}(N^\tau_0)\oplus H^*_{T_2}((\tilde N^\tau)^{T_2})$
are injective and split the Mayer-Vietoris sequences of the back-wall diagram into short exact sequences.

Let $u\in H^*_{T_2}(X^\tau)$.
We also call $u$ any of its image in Diagram \eqref{BD}, using the various
homomorphisms, including the inverses of the $\kappa$'s.
As $u=0$ in $H^*_{T_2}(\tilde N_0)$, there exists $v\in H^{2*}_T(N)$
with $v=u$ in $H^{2*}_T(N_0)\oplus H^{2*}_T(\tilde N^T)$.
Using the injectivity of
$H^*_{T_2}((X^\tau)^{T_2})\to H^*_{T_2}(N^\tau_0)\oplus H^*_{T_2}((\tilde N^\tau)^{T_2})$,
we get $u=v$ in $H^*_{T_2}((X^\tau)^{T_2})$.
By Condition~(c) and \lemref{L.CS}, there exists
$w\in H^{2*}_T(X)$ with $w=u$ in $H^*_{T_2}((X^\tau)^{T_2})$.
This proves that $A^\tau\subset \bkappafix(A)$.

To prove that $\bkappafix(A)\subset A^\tau$, let $u\in H^{2*}_T(X)$.
By a diagram chase as above,
there exists $v\in H^*_{T_2}(N^\tau)$ with $u=v$ in $H^*_{T_2}((X^\tau)^{T_2})$.
By Condition~(c) and \lemref{L.XTXT2},
$N^\tau=\skel{T_2}{1}(X^\tau)$. Using \lemref{L.CS2}, there exists
$w\in H^*_{T_2}(X^\tau)$ with $w=u$ in $H^*_{T_2}((X^\tau)^{T_2})$. This proves that
$\bkappafix(A)\subset A^\tau$.

Note that the ring homomorphism
$\kappa_T\:H^{2*}(X_T)\hfl{\approx}{} H^{*}((X_T)^\tau)$ that we have constructed
satisfies
\begin{equation}\label{qtaukap-kapq}
q^\tau\pcirc\kappa_{\scr T}=\kappafix\pcirc q \, .
\end{equation}
\end{ccote}

\begin{ccote} Construction of the ring isomorphism $\kappa\: H^{2*}(X)\to H^{*}(X^\tau)$. \rm
As $X$ is $T$-equivariantly formal, Proposition~\ref{Peqform}
tells us that the ordinary ${\rm mod\,}2$ cohomology $H^{2*}(X)$
can be recovered from the equivariant cohomology: the ring homomorphism
$\psi\:H^{2*}_{T}(X)\onto H^{2*}(X)$ descends to an isomorphism
\begin{equation}
H^{2*}_{T}(X)\otimes_{H^{2*}_{T}(pt)}\bbz_2\hfl{\approx}{} H^{2*}(X).
\end{equation}
As $X^\tau$ is $T_2$-equivariantly formal by~\eqref{Xef},
Proposition \ref{Peqform} again tells us that the ring homomorphism
$\psi^\tau\:H^{*}_{T_2}(X^\tau)\onto H^{*}(X^\tau)$ descends to
a graded ring isomorphism
\begin{equation}
H^{*}_{T_2}(X^\tau)\otimes_{H^{*}_{T_2}(pt)}\bbz_2\hfl{\approx}{} H^{*}(X^\tau).
\end{equation}
By its construction, the ring isomorphism
$\kappa_{\scr T}\: H^{2*}_{T}(X)\hfl{\approx}{} H^{*}_{T_2}(X^\tau)$
is an isomorphism of modules over the ring isomorphism
$H^{2*}_T(pt)\to H^{*}_{T_2}(pt)$. Therefore,
it descends to a graded ring isomorphism
$\kappa\:H^{2*}(X) \hfl{\approx}{}H^{*}(X^\tau)$.
With this definition, the equation
\begin{equation}\label{psitaukT-kpsi}
\psi^\tau\pcirc \kappa_{\scr T}=\kappa\pcirc \psi
\end{equation}
is satisfied.
\end{ccote}

\begin{ccote} Construction of a section
$\sigma_{\scr T}\:H^{*}(X_T)\to H^{*}_C(X_T)$ so that $(\kappa_{\scr T},\sigma_{\scr T})$ is
a \hfra\ for $(X_T,\tau)$. \ \rm
Let $(\bkappafix,\bsigmafix)$ be the \hfra\ for $X^T$.
The desired section $\sigma_{\scr T}$ will fit in the commutative diagram
$$
\begin{array}{c}{\xymatrix@C-3pt@M+2pt@R-4pt{%
H^{*}(X_T) \ar[d]^(0.50){q}
\ar@{.>}@/_/[r]_(0.50){\sigma_{\scr T}}  &
H^{*}_C(X_T) \ar[d]^(0.50){q_{\scr C}}
\ar@{>->}[r]^(0.50){r_{\scr T}} \ar[l]_(0.50){\rho_{\scr T}}  &
H^{*}_C((X_T)^\tau) \ar[d]^(0.50){q_{\scr C}^\tau}
&  H^{*}((X_T)^\tau)[u] \ar[l]_(0.50){\approx}
\ar[d]^(0.50){q^\tau[u]}
\\
H^{*}((X^T)_T) \ar@/_/[r]_(0.50){\bsigmafix}  &
H^{*}_C((X^T)_T) \ar@{>->}[r]^(0.50){\rfix} \ar[l]_(0.50){\rhofix} &
H^{*}_C((X^T)_T)^\tau) &
H^{*}(((X^T)_T)^\tau)[u]\ar[l]_(0.50){\approx}
}}\end{array}
$$
where the vertical arrows are induced by the inclusion $X^T\into X$
(the notations coincide with that of Diagram~\eqref{constrdiag}).
We have to justify that the last two
vertical arrows are injective. But, under the identifications
$$
H^*((X_T)^\tau) = H^*((X^\tau)_{T_2}) = H^*_{T_2}(X^\tau)
$$
and
$$
H^{*}(((X^T)_T)^\tau)=H^{*}(((X^T)^\tau)_{T_2})=
H^*_{T_2}((X^T)^\tau) = H^*_{T_2}((X^\tau)_{T_2})
$$
the map $ q^\tau[u]$ coincides with the homomorphism
$H^*_{T_2}(X^\tau)\to H^*_{T_2}((X^\tau)_{T_2})$ induced
by the inclusion $(X^\tau)_{T_2}\into X^\tau$.

Note that we just need to construct a section
$\sigma_{\scr T}\:H^{*}(X_T)\to H^{*}_C(X_T)$ such that
$q_{\scr C}\pcirc\sigma_{\scr T}=\bsigmafix\pcirc q$. Indeed,
if $a\in H^{2m}(X_T)$, the conjugation equation for
$(\bkappafix,\bsigmafix)$ implies
$$
q_{\scr C}^\tau\pcirc r_{\scr T}\pcirc\sigma_{\scr T}(a) =
\rfix\pcirc\bsigmafix\pcirc q(a) =\bkappafix\pcirc q(a) u^m + \lt{m} \, .
$$
As $q_{\scr C}^\tau$ is injective, this implies that
$$
r_{\scr T}\pcirc\sigma_{\scr T}(a) = \tilde a u^m + \lt{m} \, ,
$$
with $\tilde a\in H^m((X_T)^\tau)$ satisfying
$q^\tau(\tilde a)=\bkappafix(a)$. By construction of $\kappa_{\scr T}$, one has
$q^\tau\pcirc\kappa_{\scr T}(a)=\kappafix\pcirc q(a)$.
As $q^\tau$ is injective, this implies that
$\tilde a=\kappa_{\scr T}(a)$. Hence, $(\kappa_{\scr T},\sigma_T)$ satisfies the
conjugation equation and is therefore a \hfra.

As we just need to construct an additive section $\sigma_{\scr T}$, we take
the following induction hypothesis $\calh_m$: \em
for $k\leq m$, there exists a section
$\sigma_{\scr T}\:H^{2k}(X_T)\to H^{2k}_C(X_T)$ of $\rho_{\scr T}$ such that
$q_{\scr C}\pcirc\sigma_{\scr T}=\bsigmafix\pcirc q$\rm .
Hypothesis $\calh_0$ is clearly satisfied: we may assume without loss of generality that
$X$ is arc-connected; and we may then define $\sigma_T(1)=1$, where $1\in H^0(-)$ is
the unit of $H^*(-)$. Now assume by induction that $\calh_{m-1}$ holds.
The space $X_T$ is $\tau$-\ef\ by~\eqref{XTef}, so there exists a section
$\sigma_0\:H^{2m}(X_T)\to H^{2m}_C(X_T)$ of $\rho_{\scr T}$.
We have $\rhofix\pcirc q_{\scr C}\pcirc\sigma_0=q$. Therefore, for any
$a\in H^{2m}(X_T)$, we know that
$q_{\scr C}\pcirc\sigma_0(a)\equiv \bsigmafix\pcirc q(a)$ modulo $\ker\rhofix$.
This kernel is the ideal generated by $u$. As $H^{odd}_C((X^T)_T)=0$,
only even powers of $u$ occur and moreover
\begin{equation}\label{constrkT-eq30}
q_{\scr C}\pcirc\sigma_0(a) = \bsigmafix\pcirc q(a) +
\sum_{i=0}^m \sigmafix(b_{2m-2i}) u^{2i} \, ,
\end{equation}
where $b_{2j}$ are classes in $H^{2j}((X^T)_T)$ depending on the
choice of $\sigma_0$. We will modify $\sigma_0$ by
successive steps until $b_{2j}=0$ for all $j=m,m-1,\dots,0$.

The conjugation equation for $(\bkappafix,\bsigmafix)$ implies
\begin{equation}\label{constrkT-eq40}
\rfix\pcirc q_{\scr C}\pcirc\sigma_0(a) = \kappafix\pcirc q(a)u^m + \lt{m}(a) +
\sum_{i=0}^m \big(\kappafix(b_{2m-2i}) u^{m+i} + \lt{m-i}(b_{2m-2i})\big)\, .
\end{equation}
As $q_{\scr C}^\tau$ is injective, this implies that
\begin{equation}\label{constrkT-eq50}
r_{\scr T}\pcirc\sigma_0(a)= c_0 u^{2m} + \lt{m} \, ,
\end{equation}
with $c_0\in H^0((X_T)^\tau)$ satisfying $q^\tau(c_0)=\kappafix(b_0)$.
As $\kappa_{\scr T}$ is an isomorphism,
there exists $\tilde c_0\in H^0(X_T)$ with $\kappa_{\scr T}(\tilde c_0)=c_0$.
Define a new section
$$\sigma_1\:H^{2m}(X_T)\to H^{2m}_C(X_T)$$ of $\rho_{\scr T}$ by
$\sigma_1(a)=\sigma_0(a) + \sigma_{\scr T}(\tilde c_0) u^{2m}$.
By induction hypothesis,
$q_{\scr C}\pcirc\sigma_{\scr T}(\tilde c_0) = \sigmafix\pcirc q(\tilde c_0)$. By construction of $\kappa_{\scr T}$, one has
$q^\tau\pcirc\kappa_{\scr T}=\kappafix\pcirc q$.
Therefore,
\begin{eqnarray*}
\rfix\pcirc q_{\scr C}\pcirc\sigma_1(a) &=&
\rfix\pcirc q_{\scr C}\pcirc\sigma_0(a) + \rfix(q_{\scr C}\pcirc\sigma_{\scr T}(\tilde c_0))u^{2m} \\ &=&
\rfix\pcirc q_{\scr C}\pcirc\sigma_0(a)+\rfix(\sigmafix\pcirc q(\tilde c_0))u^{2m}  \\ &=&
\rfix\pcirc q_{\scr C}\pcirc\sigma_0(a)+ \kappafix\pcirc q(\tilde c_0))u^{2m}
\\ &=&
\rfix\pcirc q_{\scr C}\pcirc\sigma_0(a)+ q^\tau\pcirc\kappa_{\scr T}(\tilde c_0))u^{2m}
\\ &=&
\rfix\pcirc q_{\scr C}\pcirc\sigma_0(a)+ q^\tau(c_0)u^{2m}
\\ &=&
\rfix\pcirc q_{\scr C}\pcirc\sigma_0(a)+ \kappafix(b_0)u^{2m}
\\ &=&
\kappafix\pcirc q(a)u^m + \lt{m}(a) \\
& & + \sum_{i=0}^{m-1} \big(\kappafix(b_{2m-2i}) u^{m+i} + \lt{m-i}(b_{2m-2i})\big)\, .
\end{eqnarray*}
The injectivity of $\rfix$ implies that Equation~\eqref{constrkT-eq30}
is replaced by
\begin{equation}\label{constrkT-eq130}
q_{\scr C}\pcirc\sigma_1(a) = \bsigmafix\pcirc q(a) +
\sum_{i=0}^{m-1} \sigmafix(b_{2m-2i}) u^{2i} \, ,
\end{equation}
We thus have modified $\sigma_0$ so that $b_0=0$. Now,
using as above the injectivity of $q_{\scr C}^\tau$ this permits us
to transform~\eqref{constrkT-eq50} into
\begin{equation}\label{constrkT-eq52}
r_{\scr T}\pcirc\sigma_1(a)= c_1 u^{2m-1} + \lt{m-1} \, ,
\end{equation}
with $c_1\in H^1((X_T)^\tau)$ satisfying $q^\tau(c_1)=\kappafix(b_0)$.
Again, write $c_1=\kappa_{\scr T}(\tilde c_1)$
with $\tilde c_1\in H^2(X_T)$ and define a new section
$\sigma_2\:H^{2m}(X_T)\to H^{2m}_C(X_T)$ of $\rho_{\scr T}$ by
$\sigma_2(a)=\sigma_1(a) + \sigma_{\scr T}(\tilde c_1) u^{2m-2}$.
Proceeding as above, we prove that if we replace $\sigma_1$ by
$\sigma_2$ in~\eqref{constrkT-eq130}, the summation index runs
only till $m-2$, i.e. $b_0=b_2=0$. If we keep going this way as
long as possible, we get $\sigma_m\:H^{2m}(X_T)\to H^{2m}_C(X_T)$ with
$b_0=b_2=\cdots = b_{2m}=0$. Extending
$\sigma_{\scr T}$ in degree $2m$ by $\sigma_m$ proves $\calh_m$ holds.  So
by induction, we have our section $\sigma_T$.
\end{ccote}

\begin{ccote}\label{hfraforX} Construction of a section
$\sigma\:H^{*}(X)\to H^{*}_C(X)$ so that $(\kappa,\sigma)$ is
a \hfra. \ \rm
The relevant diagram is
$$
\begin{array}{c}{\xymatrix@C-3pt@M+2pt@R-4pt{%
H^{*}(X_T) \ar[d]^(0.50){\psi}
\ar@/_/[r]_(0.50){\sigma_{\scr T}}  &
H^{*}_C(X_T) \ar[d]^(0.50){\psi_{\scr C}}
\ar@{>->}[r]^(0.50){r_{\scr T}} \ar[l]_(0.50){\rho_{\scr T}}  &
H^{*}_C((X_T)^\tau) \ar[d]^(0.50){\psi_{\scr C}^\tau}
\\
H^{*}(X) \ar@{.>}@/_/[r]_(0.50){\sigma}
\ar@{-->}@/^/[u]^(0.50){s} &
H^{*}_C(X) \ar@{>->}[r]^(0.50){r} \ar[l]_(0.50){\rho} &
H^{*}_C(X^\tau) &
}}\end{array}
$$
Being an even-cohomology space, $X$ is $T$-\ef. We can thus choose
an additive section $s\:H^{*}(X)\to H^{*}(X_T)$ of $\psi$ and define
$\sigma\:H^{*}(X)\to H^{*}_C(X)$ by
$\sigma=\psi_{\scr C}\pcirc\sigma_{\scr T}\pcirc s$. The linear map $\sigma$
is an additive section of $\rho$ and, for $a\in H^{2m}(X)$, we have
\begin{eqnarray*}
r\pcirc\sigma(a) &=& r\pcirc\psi_{\scr C}\pcirc\sigma_T\pcirc s(a)
\\ &=&
\psi_{\scr C}^\tau\pcirc r_{\scr T}\pcirc\sigma_{\scr T}\pcirc s(a)
\\ &=&
\psi_{\scr C}^\tau\big(\kappa_{\scr T}\pcirc s(a)u^m + \lt{m}\big)
\\ &=&
\big(\psi^\tau\pcirc \kappa_{\scr T}\big)\pcirc s(a)\,u^m + \lt{m}
\\ &=&
\big(\kappa\pcirc \psi\big)\pcirc s(a)\,u^m + \lt{m}
\\ &=&  \kappa(a)u^m + \lt{m} \, .
\end{eqnarray*}
From the fourth to the fifth line, we have used that $\psi^\tau\pcirc \kappa_{\scr T}=\kappa\pcirc \psi$,
as noted in \ref{psitaukT-kpsi}. Therefore, the conjugation equation is satisfied and $(\kappa,\sigma)$ is a \hfra\ for $X$.
\end{ccote}

\noindent With \ref{hfraforX}, the proof of \thref{Prop1sk} is now complete. \cqfd

\paragraph{Proof of \corref{C.GKM}}
The hypotheses imply that the restriction of $\tau$ to an edge $E$,
which is a $2$-sphere, is conjugate to a reflection (through an equatorial plane).
We leave to the reader the details of a proof that we summarize in three steps:
(1) by an elementary argument, one shows that
$\tau$ has $2$ fixed points on each non-trivial $T$-orbit; this implies that
$E^\tau$ is a circle; (2) by the Sch\"onflies theorem, there is a homeomorphism
from $E$ to $S^2$ sending $E^\mu$ to a great circle; (3) by the Alexander trick,
the resulting involution on $S^2$ is conjugate to a reflection.
This implies that each edge is a conjugation $2$-sphere in the sense of
\cite[Example~3.6]{HHP}. Hence, each edge is a \csp\ and the hypotheses
of \thref{Prop1sk} are satisfied. \cqfd

\paragraph{Proof of \corref{C.edgeHam}}
By \cite[Remark~3.1]{HHP}, $\tau$ preserves each
arc-connected component of $X^T$. In consequence, for each edge $E$
of $X$, Hypothesis (a) of \corref{C.edgeHam} implies that $E^T$ is a \csp.
By \thref{HamilCS}, each edge is then a \csp.
The hypotheses of \thref{Prop1sk} are therefore satisfied. \cqfd

\section{Remarks}\label{S.comm}

\begin{ccote}\rm
The following example shows that the condition $\skel{T}{0}(X)=\skel{T_2}{0}(X)$
does not imply that $\skel{T}{1}(X)=\skel{T_2}{1}(X)$, even for spaces like those
occurring in \corref{C.GKM} or~\ref{C.edgeHam}.
We consider the Hamiltonian action of $S^1$ on $S^2\subset \bbc\times\bbr$ given by
$g\cdot (x,t)=(gz,t)$, compatible with the involution $(z,t)^\tau=(\bar z,t)$.
Points of $S^2$ will be denoted by $x$, $y$, etc. Let $p_\pm=(0,\pm 1)$ be the north and south poles.
Let $T=S^1\times S^1$ acting on $X=S^2\times S^2$ by
$$
(g,h)\cdot (x,y) = (gh\cdot x,gh^\mun\cdot y) \, .
$$
The fixed point sets for $T$ and $T_2$ are equal: $\skel{T}{0}(X)^\tau=\skel{T_2}{0}(X^\tau)$,
consisting of the four points $\{p_\pm\}\times\{p_\pm\}$. By \proref{Peqform3}, $X$ is
$T$-\ef\ and $X^\tau$ is $T_2$-\ef.

The $T$-equivariant $1$-skeleton is a graph of four $2$-spheres
$$
\skel{T}{0}(X)=\{(x,y)\mid x=p_\pm \hbox{ or } y=p_\pm\} \, .
$$
Therefore, $X$ is a GKM-space. But
$\skel{T}{1}(X)\neq\skel{T_2}{1}(X)$
since $\skel{T_2}{1}(X)=X$.
Also, $\skel{T}{1}(X)^\tau\neq\skel{T_2}{1}(X^\tau)$
since $\skel{T_2}{1}(X^\tau)=X^\tau$.
\end{ccote}

\begin{ccote}\rm The condition $\skel{T}{i}(X)=\skel{T_2}{i}(X)$ for $i=0,1$ of our main
theorems is already implicitly present in earlier papers \cite{Sd, BGH}
which are dealing with GKM Hamiltonian manifolds. In \cite{Sd}, one requires that for each
point of $x\in X^T$, the characters involved in the $2$-spheres adjacent to $x$
are pairwise independent over $\bbz_2$.
In \cite[, p.~373]{BGH}, one asks that $X^T=X^{T_2}$ and that ``the real locus of the one-skeleton
is the same as the one-skeleton of the real locus''.
In general, these conditions
are weaker than $\skel{T}{i}(X)=\skel{T_2}{i}(X)$ for $i=0,1$ (see \lemref{L.XTXT2})
but, they are equivalent for a GKM Hamiltonian manifold.
To see this, work with the local normal around a $T$-fixed point; in this model
the $T$-action and the involution are linear.

\end{ccote}

\sk{3}\noindent{\bf Acknowledgments: }\rm The authors thank Volker Puppe for useful discussions,
in particular about Lemmas~\ref{L.CS2} and~\ref{L.CS}.

\end{document}